# Three Theorems on Modular Sieves that suggest the Prime Difference is

$$O(\pi(p(n)^{1/2}))$$

Bhupinder Singh Anand[1]

Although this 1964 paper is not, strictly, a foundational paper, it developed as an off-shoot to the foundational query: Do we discover the natural numbers (Platonically), or do we construct them linguistically? The paper also assumes computational significance in the light of Agrawal, Kayal and Saxena's August 2002 paper, "PRIMES is in P", since both the Trim and Compact Number Generating algorithms - each of which generates all the primes - are deterministic algorithms that run in polynomial time, and suggest that the Prime Difference, $d_P(n)$, is $O(\pi(p(n)^{1/2}))$.

## Contents

[1] The author is an independent scholar. E-mail: re@alixcomsi.com; anandb@vsnl.com. Postal address: 32, Agarwal House, D Road, Churchgate, Mumbai - 400 020, INDIA. Tel: +91 (22) 2281 3353. Fax: +91 (22) 2209 5091.

## 1: A PRIME NUMBER GENERATING THEOREM (1964)

We address the question:

Do we necessarily discover the primes, as in the sieve of Eratosthenes, or can we also generate them sequentially?

In other words, given the first $k$ primes, can we generate the prime $p(k+1)$ algorithmically, without testing any number larger than $p(k)$ for primality?

We give an affirmative answer, defining two algorithms based on the following theorem (where $i$, $j$, $k$, … represent natural numbers):

**Theorem 1.1**: For any given natural number $k$, and all $i < k$, let $p(i)$ denote the $i$'th prime, and define the sequence $\{a_P(i): i = 1 \text{ to } k\}$ such that $a_P(i) < p(i)$, and:

$$p(k) + a_P(i) \equiv 0 \ (\text{mod } p(i))$$

Let $d_P(k)$ be such that, for all $l < d_P(k)$, there is some $i$ such that:

$$l \equiv a_P(i) \ (\text{mod } p(i)).$$

Then $d_P(k)$ is the Prime Difference, defined by:

$$p(k+1) - p(k) = d_P(k).$$

**Proof**: For all $l < d_P(k)$, there is some $i$ such that:

$$p(k) + l \equiv 0 \ (\text{mod } p(i)).$$

Since, by Bertrand's Postulate ([HW60], p343), $p(k+1) < 2p(k)$ for all natural numbers $k$, it follows that:

$$d_P(k) < p(k)$$

Hence, $p(k) + d_P(k) < p(k)^2$, and so it is the prime $p(k+1)$□

## 2: TRIM NUMBERS (1964)

The significance of the Prime Number Generating Theorem is seen in the following algorithm.

We define Trim numbers recursively by $t(1) = 2$, and $t(n+1) = t(n) + d_T(n)$, where $p(i)$ is the $i$'th prime and:

(i)    $d_T(1) = 1$, and $a_T(2, 1) = 1$:

(ii)    $d_T(n)$ is the smallest even integer that does not occur in the $n$'th sequence:

$$\{ a_T(n, 1), \ldots, a_T(n, n\text{-}1)\};$$

(iii)    $j_i \geq 0$ is selected so that, for all $0 < i \leq (n\text{-}1)$:

$$0 \leq a_T(n+1, i) = \{( a_T(n, i) - d_T(n) + j_i * p(i)\} < p(i)$$

It follows that the Trim number $t(n+1)$ is, thus, a prime unless all its prime divisors are less than $d_T(n)$.

(In [An64] I illustrate how Trim numbers are generated sequentially as a 2-dimensional modular sieve.)

## 2.1: A TRIM NUMBER THEOREM … (1964)

**Theorem 2.1:** For all $n > 1$, $t(n) < n^2$.

**Proof:** $t(n) - \{t(n\text{-}1) = d_T(n\text{-}1)\}$

$d_T(n\text{-}1) < 2(n\text{-}1)$

$t(n) - t(n\text{-}1) < 2(n\text{-}1)$

$t(n) - 2 < 2(n\text{-}1) + 2(n\text{-}2) + \ldots + 2(1)$

$t(n) < [2 + 2\{(n\text{-}1)*n/2\}]$

$t(n) < (n^2 - n + 2)$

Hence, for all $n > 1$, $t(n) < n^2 \square$

## 2.2: ... AND TWO CONJECTURES (1964)

**Conjecture 1.1:** For all $k$, $t(k) < p(n) < t(k+1)$ for some $n$.

**Conjecture 1.2:** If $p(n+1) = \{p(n) + d_p(n)\}$, then $d_p(n) = O(n)$.

## 2.3: IS THE THEOREM SIGNIFICANT? (2005)

The significance of the above theorem lies in the following observation [Ca05]:

"Is there always a prime between $n^2$ and $(n+1)^2$? In 1882 Opperman stated $\pi(n^2+n) > \pi(n^2) > \pi(n^2-n)$, $(n>1)$, which also seems very likely, but remains unproven ([Ri95], p248).

Both of these conjectures would follow if we could prove the conjecture that the prime gap following a prime $p$ is bounded by a constant times $(\log p)^2$."

## 3: COMPACT NUMBERS (1964)

Compact numbers are defined recursively by $c(1) = 2$, and $c(n+1) = c(n) + d_c(n)$, where $p(i)$ is the $i$'th prime and:

(*i*)   $d_c(1) = 1$, and $a_c(2, 1) = 1$;

(*ii*)   $d_c(n)$ is the smallest even integer that does not occur in the $n$'th sequence:

$$\{a_c(n, 1), \ldots, a_c(n, k)\};$$

(*iii*)   $j_i \geq 0$ is selected so that, for all $0 < i \leq k$:

$$0 \leq a_c(n+1, i) = \{ a_c(n, i) - d_c(n) + j_i*p(i) \} < p(i);$$

(*iv*)   $k$ is selected so that:

$$p(k)^2 < c(n) \leq p(k+1)^2;$$

(*v*)   if $c(n) = p(k+1)^2$, then:

$$a_c(n, k+1) = 0.$$

It follows that the compact number $c(n+1)$ is either a prime, or a prime square, unless all, except a maximum of 1, prime divisors of the number are less than $d_c(n)$.

(In [An64] I illustrate how Compact numbers are generated sequentially as a significantly more compact, 2-dimensional, modular sieve.)

### 3.1: TWO COMPACT NUMBER THEOREMS … (1964)

**Theorem 3.1:** There is always a compact number $c(m)$ such that:

$$n^2 < c(m) \leq (n+1)^2.$$

**Proof:** If $p(k) \leq n < p(k+1)$, then $(n+1) \leq p(k+1)$.

Now $\{(n+1)^2 - n^2\} = \{2n + 1\} > 2p(k) > 2k$.

Hence $(n+1)^2 > \{n^2 + 2k\}$.

Since there is always a compact number in any interval larger than $2k$ between $p(k)^2$ and $p(k+1)^2$, the theorem follows$\square$

**Theorem 2:** For sufficiently large $n$, $d_c(n) <$ constant*$\{c(n)/\log c(n)\}^{1/2}$

**Proof:** By definition, $d_c(n) \leq 2k$, where $p(k)^2 < c(n) \leq p(k+1)^2$.

The theorem follows from the Prime Number Theorem, since $k$ is the number of primes less than $c(n)^{1/2}$$\square$

### 3.2: ... AND TWO MORE CONJECTURES (1964)

**Conjecture 3.1:** If $p(k)^2 < p(n) < p(k+1)^2$, then:

$$\{p(n+1) - p(n)\} = O(k)$$
$$= O(\pi(p(n)^{1/2})).$$

**Conjecture 3.2:**

$$\lim_{n \to \infty}^{\sup} \{\text{Number of primes} < n\}/\{\text{Number of compacts} < n\} = 1/2.$$

*(Transcribed from original notes: Thursday 8th December 2005 2:17:06 AM by re@alixcomsi.com. Updated: Saturday 24th December 11:01:33 AM by re@alixcomsi.com.)*